\documentclass[11pt,conference]{ukacm2024}
\usepackage{floatrow}
\newfloatcommand{capbtabbox}{table}[][\FBwidth]

\usepackage{blindtext}
\usepackage{parskip}
\usepackage{pslatex}
\usepackage{graphicx}
\usepackage{amsmath}
\usepackage{amsfonts}
\usepackage{amssymb}
\usepackage{hyperref}
\usepackage[T1]{fontenc}
\usepackage[utf8]{inputenc}
\usepackage{babel}
\usepackage[font=small,labelfont=bf]{caption}

\title{Study of Landslides through a Stabilised Semi-Implicit Material Point Method}

\author{M. Xie$^{1*}$, P. Navas$^2$, S. L\'opez-Querol$^1$}

\heading{M. Xie, P. Navas and S. L\'opez-Querol}

\address{$^{1}$ Dept. of Civil, Environmental \& Geomatic Engineering, University College London, UK.\\
mian.xie.18@ucl.ac.uk, s.lopez-querol@ucl.ac.uk \\
$^{2}$ Dept. Continuum Mechanics and Theory of Structures, UPM, Madrid, Spain.
pedro.navas@upm.es}

\abstract{In this research, a new semi-implicit two-phase double-point material point method is proposed, in which the soil and water phases are modelled using two distinct sets of material points, both being stabilised with a novel approach. The Nor-Sand constitutive model is implemented to simulate more realistic soil behaviour. Some landslide numerical examples are presented to investigate the performance of the proposed method and highlight the importance of using the double-point approach. The formulation with two sets of material points shows significantly different but more reliable results in the cases of landslides, compared with the conventional single-point approach. Furthermore, this research shows that the additional computational cost given by the additional water material points is acceptable. Therefore, it is recommended to use two sets of material points for some large deformation geotechnical problems.}


\begin{document}

\section{Introduction}

The Material Point Method (MPM) has become a popular continuum method for modelling large deform-ation geotechnical problems such as landslides. The MPM formulation is very similar to that of the Finite Element Method (FEM), except that the iteration points (i.e., material points) can move independently from the mesh,  which is set to its original position at every step. Special treatment of iteration points and mesh allows MPM to model large deformation problems without mesh distortion; however, instabilities, such as quadrature error, cell-crossing noise and volumetric locking, arise. Although many researchers have attempted to solve these instabilities, fully stabilising the MPM is still a challenging task. Due to its unstable nature, MPM studies usually adopt explicit approaches and basic constitutive models.

To ensure convergence in an explicit approach, the time step must be smaller than the critical time step controlled by the Courant–Friedrichs–Lewy (CFL) condition. The critical time step defined by this condition is related to the bulk modulus of the material. Pore water is usually an important factor in geotechnical problems. However, the bulk modulus of water is about 100 times larger than that of the soil, resulting in an extremely small critical time step in a traditional explicit soil-water coupled MPM formulation. The time step of an implicit approach is not limited by the CFL condition, but the convergence issue can easily be raised when using an advanced constitutive model. Recently, the semi-implicit two-phase MPM has been derived based on the incremental fractional step method \cite{Kularathna2021,Yuan2023}. In this approach, the water phase is solved implicitly, resulting in a significantly larger time step compared to the fully explicit approach. However, the available methods \cite{Kularathna2021,Yuan2023} use a single set of material points to model both the soil and water phases (so-called single-point approach) with a very basic constitutive model. We have derived a new semi-implicit two-phase double-point MPM, in which the soil and water phases are modelled using two distinct sets of material points. The proposed method is stabilised by a newly derived modified F-bar method aiming to stabilise advanced constitutive models, in which both deformation gradient increment and last converged elastic left Cauchy-Green tensors are stabilised. The Nor-Sand constitutive model is implemented to simulate the soil behaviour in a more realistic way. The derivation and formulation are presented in our extended paper \cite{Xie2024}. 

In this research, some landslide examples are presented, to investigate the performance of the modified F-bar method and highlight the importance of using the double-point approach. More numerical examples, validations and the method of B-spline MPM can be found in the extended paper \cite{Xie2024}.

\section{Numerical examples}

\subsection {Landslides - dry condition}
\begin{minipage}{\textwidth}
    \begin{minipage}[b]{0.49\textwidth}
        \centering
        \includegraphics[width=0.93\textwidth]{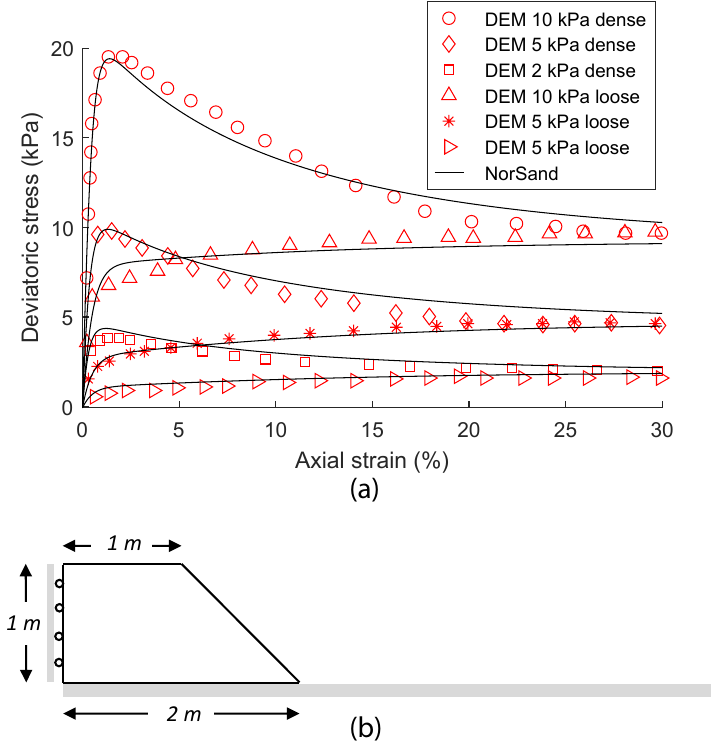}
        \captionof{figure}{(\textbf{a}) Calibration of Nor-Sand with DEM drained triaxial; (\textbf{b}) Illustration of numerical model \cite{Xie2024}.}\label{fig1}
    \end{minipage}
    \hfill
    \begin{minipage}[b]{0.49\textwidth}
        \centering
        \begin{tabular}{@{}ll@{}}
            \hline
            Parameter & Value \\
            \hline
		  Soil density $\rho_s [kg/m^3]$    & 2,500 \\
		  Shear modulus $G[kPa]$    & 4,000 \\
		  Swelling index $\kappa$    & 0.002 \\
            Reference specific volume $v_{c0}$    & 1.88 \\
            Compression index $\lambda$    & 0.012 \\
		  Slope of the critical state line $M$    & 0.7097 \\
		  Yield function constant $N$    & 0.3 \\
            Plastic potential constant $\Bar{N}$    & 0.3 \\
		  Hardening coefficient $h$    & 300 \\
            Maximum dilatancy multiplier $\alpha$    & -4.5 \\
            Initial specific volume $v_0$ (dense)    & 1.629 \\
            Initial specific volume $v_0$ (loose)   & 1.736 \\
            Gravitational acceleration $g[m/s^2]$    & -9.81 \\
            \hline
            %Water density $\rho_w [kg/m^3]$    & 1000 \\
		  Initial porosity $n_0$    & 0.4 \\
		  Initial hydraulic conductivity $k_0 [m/s]$    &     $5\times10^{-3}$ \\
            \hline
        \end{tabular}
        \captionof{table}{Parameters for Landslides with Nor-Sand model.}\label{tab1}
    \end{minipage}
\end{minipage}

Lu et al. \cite{Lu2022} studied the problems of landslides using the discrete element method (DEM). DEM simulations \cite{Lu2022} were conducted using dry granular materials with different particle shapes. Also, drained triaxial tests were conducted for these particles using DEM. %Lu et al. \cite{Lu2022} have also conducted drained triaxial tests for these particles using DEM. 
In this research, we focus on spherical particles without rolling resistance \cite{Lu2022}. The parameters of the Nor-Sand constitutive model are calibrated with the DEM drained triaxial tests using the single-element finite strain driver derived by Xie et al. \cite{Xie2024}. Both  Nor-Sand and the DEM stress-strain curves %given by the Nor-Sand model and the DEM simulations 
are presented in Fig. \ref{fig1}(a). The calibrated Nor-Sand parameters for the landslides simulations are summarised in Tab. \ref{tab1}.

Fig. \ref{fig1}(b) shows the geometry and boundary condition of the numerical model. The roller boundary condition is applied on the sides of the background mesh, and the fixed boundary condition is applied on the bottom of the background mesh to represent a rough surface. A 0.05 m wide square mesh with $4^{2}$ material point per cell is used. The mesh is further refined to 0.025 m for the sensitive analysis purpose. The single-phase B-spline MPM is used with the Nor-Sand constitutive model. The total simulation time is 6 s for both dense and loose conditions, and the landslides are almost static at this time.   

\begin{figure}[t]
\centering
\includegraphics[width=0.96\textwidth]{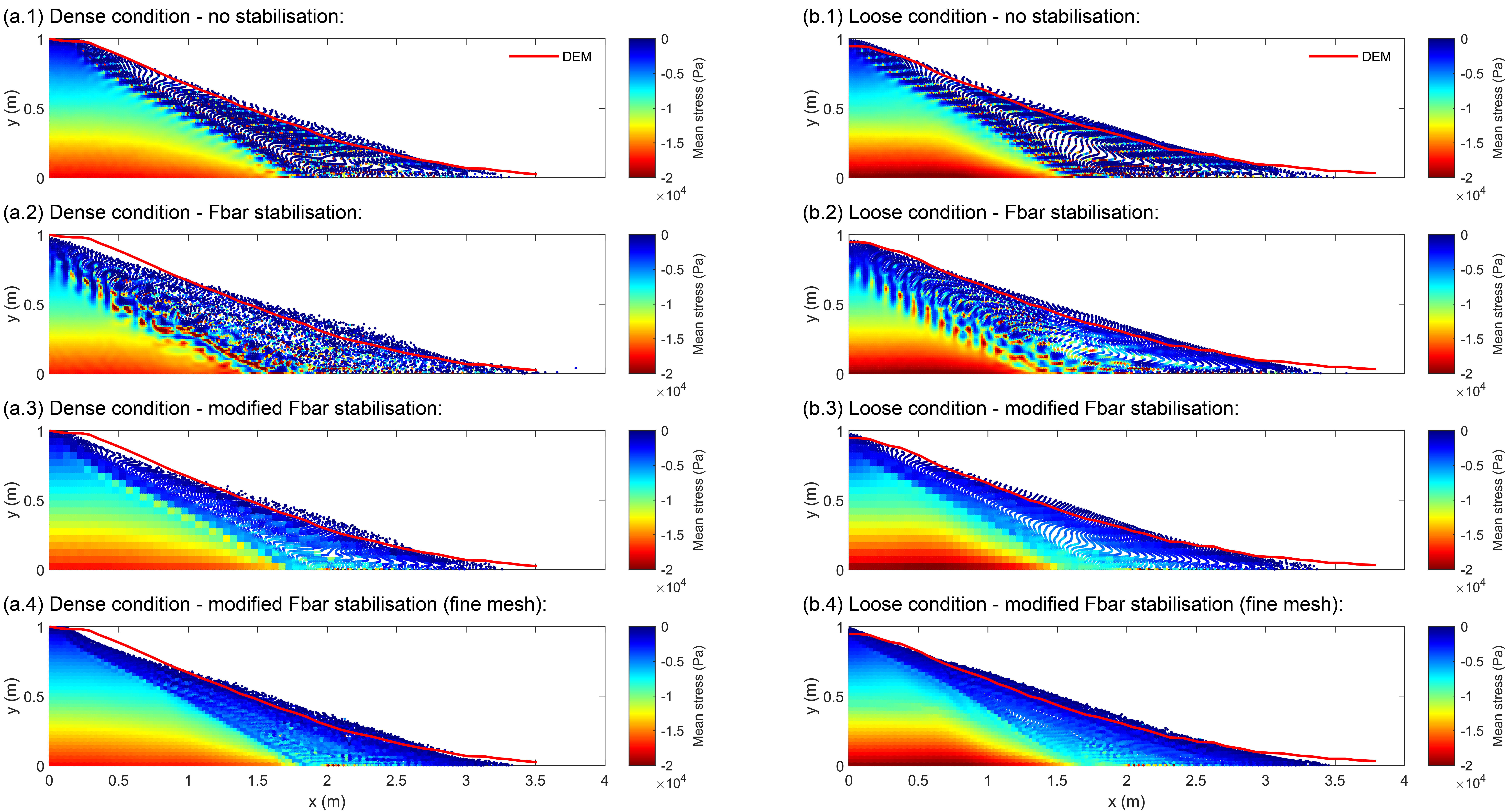}
\caption{Comparison of different stabilisation methods under the (\textbf{a}) dense condition (\textbf{b}) loose condition}\label{fig2}
\end{figure}

Fig. \ref{fig2} presents the simulation results of the proposed method with different stabilisation methods. As we can see, the standard F-bar method is not able to fully stabilise the stress oscillations under very large deformation. The performance of the standard F-bar method is even worse in the case of dilative granular material (i.e., dense condition), as shown in Figs. \ref{fig2}(a.2) and (b.2). On the contrary, the modified F-bar approach proposed by Xie et al. \cite{Xie2024} can sufficiently eliminate stress oscillations. A higher resolution stress contour can be obtained by refining the mesh as shown in Figs. \ref{fig2}(a.4) and (b.4). However, a 0.05 m mesh is sufficient to ensure the accuracy of this numerical example. In addition, Fig. \ref{fig2} shows an excellent agreement between DEM and B-spline MPM with the calibrated Nor-Sand constitutive model.

% In the modified F-bar method, both deformation gradient increment and last converged elastic left Cauchy-Green tensors are stabilised.

\begin{figure}[t]
\centering
\includegraphics[width=0.96\textwidth]{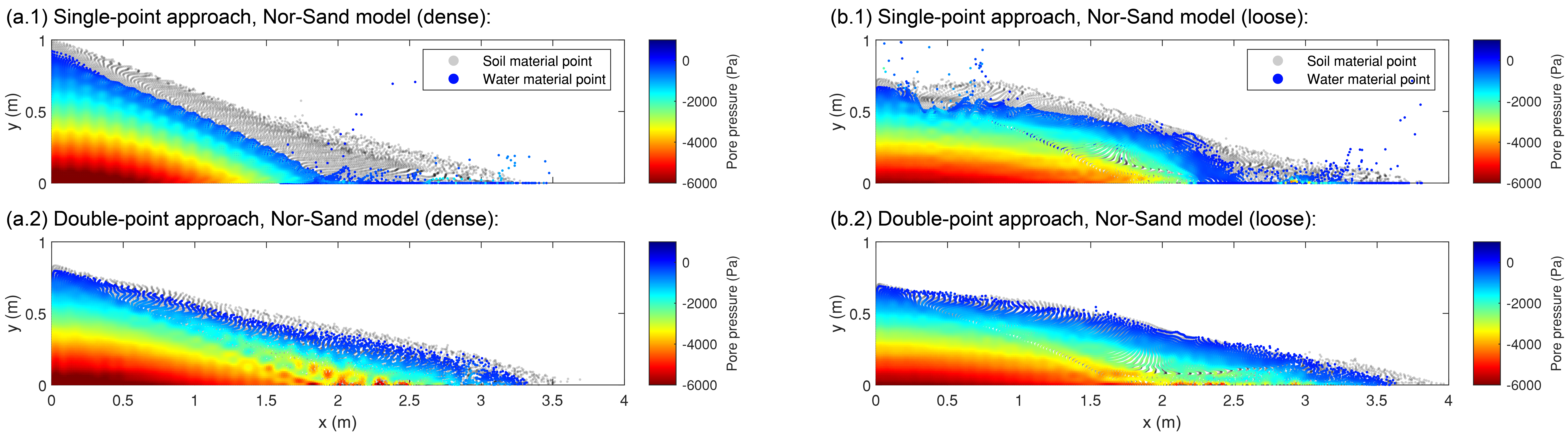}
\caption{Comparison of the single- and double-point methods under the (\textbf{a}) dense condition (\textbf{b}) loose condition.}\label{fig3}
\end{figure}

\subsection {Landslides - fully saturated condition}

For landslides with fully saturated condition, the same geometry and boundary conditions are applied as in the previous dry condition. A 0.05 m wide square mesh with $4^{2}$ material point per cell is used. In addition, the same Nor-Sand parameters are used. The additional parameters for two-phase MPM are documented in Tab. \ref{tab1}. The model initialisation process follows Xie et al. \cite{Xie2024}. The total time, excluding the initialisation, is 12 s, when the propagation of the landslide becomes static in all cases. 

The simulations are conducted using both semi-implicit single- and double-point approaches stabilised by the modified F-bar method. The soil's effective stress contour is omitted for these fully saturated landslides because it follows a similar pattern as that in the previous dry cases. An interested reader may refer to Xie et al. \cite{Xie2024} where the soil stress contour is plotted. The pore water pressure contours are presented along with the soil material points. As shown in Fig. \ref{fig3}, the soil material points are represented as grey dots underneath the water material points. Although water material points are attached to soil material points in a single-point approach, we can still graphically present the water and soil material points separately. This is because in the water material points have their own velocities, and their position can be updated based on the velocities and initial coordinates.

As we can see in Fig. \ref{fig3}, the water phases of the single-point approach experience extensive instability in the case of the Nor-Sand model under both dense and loose conditions. On the contrary, the double-point approach is very stable. This phenomenon has not been reported in previous single-point MPM studies and may be due to two reasons: (a) a basic constitutive model (e.g. Mohr-Coulomb and Drucker-Prager) is used that cannot capture a realistic soil behaviour and (b) the pore pressure contour is plotted on the soil material points. Despite instabilities, the single-phase approach shows significantly stiffer behaviour compared to the double-point approach for the landslide problem, as shown in Fig. \ref{fig3}. The landslides are not fully generated in a single-phase approach with the Nor-Sand model, resulting in an underestimated runout distance and an overestimated angle of repose. In addition to accuracy, the computational efficiency of the double-point approach has also been assessed. Single- and double-point simulations have been repeated in the same computational environment. The additional computational cost of the double-point approach fluctuates from 14\% to 16\%, which is acceptable.

\section{Conclusions}
%This research studies landslides under dry and fully saturated conditions using a stabilised semi-implicit two-phase double-point MPM with a Nor-Sand constitutive model. Using the single-point approach is not conservative for landslides problems. 
The modified F-bar method shows promising performance on MPM with a Nor-Sand model under very large deformation, whereas the F-bar method fails to stabilise the problem in this scenario. The MPM with Nor-Sand model can reproduce the DEM simulations for landslide problems as long as the material parameters are carefully calibrated. The two-phase double-point MPM is more stable and accurate than the single-point approach. In addition, the double-point approach only increases the computational cost by about 15\%.

\end{document}